%% file: korea.tex
%%<pre>
%% neumann.tex: Paper for the Korea Conference Proceedings
%% This is a Plain TeX file. 
%% Needed files:  
%% degmacs.tex:  The de Gruyter Macro file
%% neumfig.eps:  A PostScript figure for the paper
%% amssym.def:  The ams font definition file
%% epsf.tex:  To include the Postscript figure
%%    This assumes dvips is being used as the PostScript printer driver.
\input degmacs.tex
\input amssym.def
\input epsf.tex
%% Corrections to plain.tex insert macros stripped from newinsert.tex 
%% on the CTAN archives (plain.tex has bugs).
\topskip=\baselineskip
\chardef\newinsCatAt\the\catcode `\@
\catcode `\@=11
\newskip\insertskipamount\newskip\inserthardskipamount
\insertskipamount 6pt plus2pt \inserthardskipamount 6pt
\def\insertskip{\vskip\insertskipamount}
\newskip\LastSkip\def\SaveLastSkip{\LastSkip\lastskip}
\def\RestoreLastSkip{\nobreak\vskip-\LastSkip\vskip\LastSkip}
\newcount\SplitTest%        will be set to -1 if a topinsert has split
\def\SetSplitTest{\SplitTest\insertpenalties
  \insert\topins{\floatingpenalty1}%
  \advance\SplitTest-\insertpenalties}
\def\midinsert{\par
 \SaveLastSkip\penalty-150\SetSplitTest\RestoreLastSkip
 \ifnum\SplitTest=-1
  \@midfalse\p@gefalse\else\@midtrue\fi\@ins}
\def\@ins{\par\begingroup\setbox\z@\vbox\bgroup%
  \vglue\inserthardskipamount}
\def\endinsert{\egroup % finish the \vbox
  \if@mid \dimen@\ht\z@ \advance\dimen@\dp\z@
    \advance\dimen@\insertskipamount%
    \advance\dimen@\pagetotal\advance\dimen@-\pageshrink
    \ifdim\dimen@>\pagegoal\@midfalse\p@gefalse\fi\fi
  \if@mid%
    \ifdim\lastskip<\insertskipamount\removelastskip\insertskip\fi
    \nointerlineskip\box\z@\penalty-200\insertskip
  \else%
    \SaveLastSkip%                         
    \insert\topins{\penalty100 % floating insertion
    \splittopskip\z@skip
    \splitmaxdepth\maxdimen \floatingpenalty\z@
    \ifp@ge \dimen@\dp\z@
    \vbox to\vsize{\unvbox\z@\kern-\dimen@}% depth is zero
    \else \box\z@\nobreak\insertskip\fi}% 
    \RestoreLastSkip% 
   \fi\endgroup}
\catcode `\@=\newinsCatAt
%% end newinsert.tex

%% Local definitions:
\let\mop\mathop
\def\Isom{\mop{Isom}}
\def\PSL{\mop{PSL}}
\def\H{\Bbb H}
\def\C{\Bbb C}
\def\cite#1{[{\bf#1}]}
\def\author{Walter D. Neumann}
\def\title{Kleinian Groups Generated by Rotations}

\runningheads\title\author

\Title\title

\Author\author

\Abstract We discuss which Kleinian groups are commensurable with
Kleinian groups generated by rotations, with particular emphasis on
Kleinian groups that arise from Dehn surgery on a knot.

\Section Introduction

In a problem session, organized by A. Kim, of the 1991
German-Korean-SEAMS Conference on Geometry, E. Vinberg asked for a
cocompact Kleinian group which is not commensurable with a group
generated by rotations (elements of finite order). Examples are, in
fact, not hard to find. In this note we describe in some detail which
Kleinian groups have this property among the Kleinian groups that
occur as fundamental groups of Dehn surgeries on knots.  We also
briefly discuss some related questions.  \iffalse This paper is an
improved version of a report that has apparantly appeared in the
Proceedings of the above conference, published in 1994 by
Chulalungkorn University Press, but I have not yet been able to obtain
a copy.\fi

In the following $\Gamma$ and $\Lambda$ will always denote Kleinian
groups of finite covolume, that is, discrete subgroups of
$\PSL(2,\C)=\Isom^+(\H^3)$ such that the orbifold $\H^3/\Gamma$
or $\H^3/\Lambda$ has finite volume.  They are {\it commensurable} if
they have isomorphic subgroups of finite index.  By Mostow-Prasad
rigidity this is equivalent to the condition that they can be
conjugated within $\PSL(2,\C)$ so their intersection has finite
index in each.

Kleinian groups commensurable with groups generated by reflections
(rather than rotations) have been studied by E. Vinberg \cite{V} and
E.~M.~Andreev \cite{An}.  They are a very restricted class of Kleinian
groups. For example, as is pointed out in \cite{NR}, the invariant
trace field of such a Kleinian group --- indeed, of any Kleinian group
commensurable with a non-orientation-preserving subgroup of
$\Isom(\H^3)$ --- has to be preserved by complex conjugation (and the
same for the invariant quaternion algebra), which is a rare occurrence.  
For a Kleinian group
generated by rotations I know of no similar restriction on the
invariant trace field, although the invariant trace field severely
restricts the possible orders of the rotations involved (cf.~Theorem 3
below).

In the following section we discuss Dehn surgery on knots and state
and prove our main theorem (Theorem 1). In Section 2 we discuss the
case of knot complements themselves: a hyperbolic knot complement is
commensurable with a $\H^3/\Lambda$ with $\Lambda$ generated by
rotations if the knot is invertible and conjecturally in just one
other case (Theorem 2). In Section 3 we make some additional comments
about what rotations can be contained in Kleinian groups commensurable
with a given Kleinian group $\Gamma$.  If such rotations do not occur
in $\Gamma$ itself, it is reasonable to call them ``hidden rotations''
of $\Gamma$.  We show that ``most'' Kleinian groups have no hidden
rotations of orders other than $2,3,4$, or $6$.

\Section 1.~~Dehn surgery on knots

Let $(S^3,K)$ be a hyperbolic knot (i.e., $M=S^3-K$ admits a complete
finite-volume hyperbolic structure).  Thurston's hyperbolic Dehn
surgery theorem (\cite{T}, \cite{N-Z}) implies that by excluding
finitely many integer pairs $(p,q)$ we may assume the result $M(p,q)$
of $(p,q)$-Dehn-Surgery on $(S^3,K)$ is hyperbolic.  Say
$M(p,q)=\H^3/\Gamma(p,q)$.  Group-theoretically, $\Gamma(p,q)$ is the
result of factoring $\pi_1(S^3-K)$ by the normal closure of the
element $m^pl^q$, where $m,l\in\pi_1(S^3-K)$ are represented by a meridian
and a longitude in a torus neighbourhood boundary of $K$ in $S^3$.

\Proclamation Theorem 1.
\\i) If $(S^3,K)$ is non-invertible and does not branched cyclic
cover a torus knot, then for all but finitely many values of $p/q$ the
above group $\Gamma(p,q)$ is not commensurable with a group generated
by rotations.
\\ii) If $(S^3,K)$ is an invertible knot, then each $\Gamma(p,q)$
has index $2$ in a group generated by rotations.
\\iii) If $(S^3,K)$ is a $d$-fold cyclic cover of the unknot, then
each $\Gamma(p,q)$ is commensurable with a group generated by
rotations.
\\iv) If $(S^3,K)$ is not in one of the above cases, then
infinitely many $\Gamma(p,q)$ are commensurable with groups generated
by rotations and infinitely many are not.

\noindent{\bf Example}.~~
Figure 1 shows a knot $3$-fold cyclic covering an unknot.  This knot
(and its obvious generalizations) provides an example for case (iii).
However, in some sense ``most'' knots have no symmetries at all, so
they satisfy case (i).

\midinsert\centerline{\epsfxsize=\hsize\epsffile{neumfig.eps}}
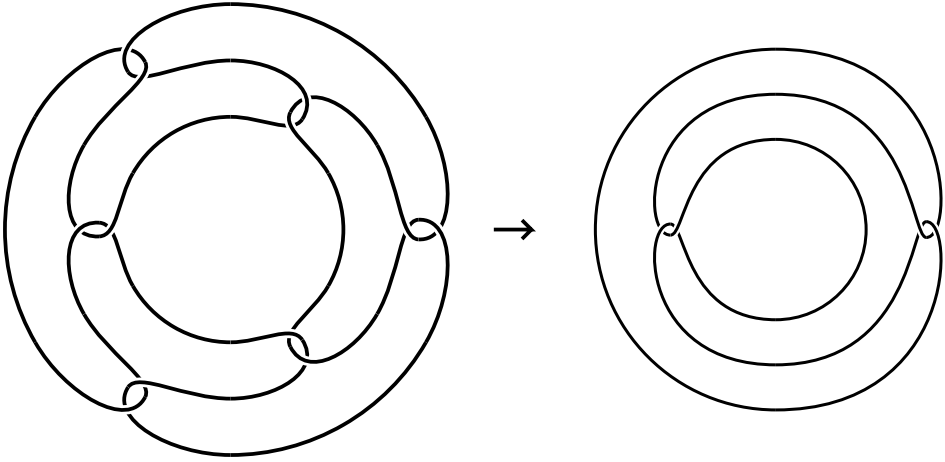
\medskip\centerline{Figure 1}
\endinsert

\Proof
It is well known (cf.~\cite{Ar}) that $\Lambda\subset\Isom^+(\H^3)$
is generated by rotations if and only if the space $\H^3/\Lambda$ is
simply connected.

Let $M(p,q)=\H^3/\Gamma(p,q)$ be as in the theorem. We shall
abbreviate $\Gamma=\Gamma(p,q)$. We first assume $(S^3,K)$ is
non-invertible.

Thurston's Dehn surgery theorem (loc.~cit.) also says ${\rm
vol}(M(p,q))<{\rm vol}(M)$, and since, by Borel \cite{B}, only
finitely many hyperbolic orbifolds with volume below a given bound are
arithmetic, we may, by excluding finitely many $(p,q)$, assume
$M(p,q)$ is non-arithmetic.  Then by Margulis (cf. \cite{Z}, ch. 6)
there exists a maximal element $\Gamma_0$ in the commensurability
class of $\Gamma$.  Denote $M_0(p,q)=\H^3/\Gamma_0$.  The map
$M(p,q)\to M_0(p,q)$ is a covering map of orbifolds. We shall show
first that for $(p,q)$ sufficiently large it is a cyclic covering.

There is a bound on the degree of $M(p,q)\to M_0(p,q)$, independent of
$(p,q)$, since a complete hyperbolic orbifold has volume above some
fixed positive bound (Margulis).  If $(p,q)$ is sufficiently large,
then the geodesic $\gamma$ added to $M$ by Dehn surgery is much
shorter than the other closed geodesics of $M(p,q)$.  It follows that,
for $(p,q)$ sufficiently large, the image $\gamma_0$ of $\gamma$ in
$M_0(p,q)$ is the shortest closed geodesic of $M_0(p,q)$, and that no
other geodesic of $M(p,q)$ could cover $\gamma_0$.  Then
$M=M(p,q)-\gamma$ covers $M_0(p,q)-\gamma_0$.  We exclude the finitely
many $(p,q)$ for which this is not the case.

Now, consider this covering restricted to the boundary $T$ of a solid
torus neighborhood of $\gamma$.  The image $T_0$ in
$M_0(p,q)-\gamma_0$ is an orbifold covered by the torus.  Let $A$ and
$A_0$ be the fundamental groups of $T$ and $T_0$.  If $A$ is not
normal in $A_0$, that is the covering $T\to T_0$ is not regular, then
$T_0$ must be a triangle orbifold of type $(2,4,4)$, $(2,3,6)$, or
$(3,3,3)$, since these are the only orbifolds covered by a torus for
which the orbifold fundamental group $A_0$ contains non-normal
torsion-free subgroups.  However, such a covering cannot extend to a
solid torus and the covering $T\to T_0$ extends to a tubular
neighborhood of $\gamma$.  Hence the covering $T\to T_0$ is regular.
As described for instance in \cite{R}, it follows that the covering
$M=M(p,q)-\gamma\to M_0(p,q)-\gamma_0$ is regular, since $\pi_1(T)$
normally generates $\pi_1(M)$ ($M$ is a knot complement); we repeat
the argument for completeness.

We show $\Gamma$ is normal in $\Gamma_0$ by showing that it is
the normal closure of a suitable subgroup.  If we identify $A_0$ with
its image in $\Gamma_0$ then $A=A_0\cap\Gamma$. The degree of our
covering is on the one hand equal to $|A_0/A|=|A_0/(A_0\cap\Gamma)|=
|A_0\Gamma/\Gamma|$ and on the other hand equal to
$|\Gamma_0/\Gamma|$, so $A_0\Gamma=\Gamma_0$.  The normal closure of
$A$ in $\Gamma_0=A_0\Gamma$ therefore equals the normal closure of $A$
in $\Gamma$, which is $\Gamma$ itself, as pointed out above.

The covering transformation group $G$ for $M\to M_0(p,q)-\gamma_0$ is
$G=\Gamma_0/\Gamma=A_0/A$, so it is the same as for $T\to T_0$.  Since
the longitude of the knot complement $M(p,q)-\gamma =S^3-K=M$
generates the kernel of $H_1(T)\to H_1(M)$, it is preserved by the
$G$-action up to sign.  Hence the same is true for the meridian, so
the $G$-action extends to an action on $S^3$.  By the solution to the
Smith Conjecture (cf.~\cite{MB}) $G$ is cyclic or dihedral.  But dihedral
is excluded by our assumption that $K\subset S^3$ is a non-invertible
knot, so $G$ is cyclic.

Now suppose that $\Lambda\subset\Gamma_0$ exists so that $\Lambda$ is
generated by rotations, i.e., the space $\H^3/\Lambda$ is simply
connected (cf.\ first sentence of this section).  Then
$\H^3/\Gamma_0$ must have finite fundamental group.  But
$\H^3/\Gamma_0$ is the result of a Dehn filling of $M/G$, where $G$
is the above cyclic group.  By \cite{BH}, since we assumed $M/G$ is
not a torus knot complement, such a Dehn surgery can give finite
fundamental group for at most $24$ values of $p/q$.  Thus by excluding
these $(p,q)$ we can avoid this, so part (i) of the theorem is proved.

(We note that this final exclusion may be necessary.  For example, if
$(S^3,K)$ is the $(-2,3,7)$-pretzel knot then $(r,1)$-Dehn surgery
gives a manifold with finite fundamental group for $r=17,18,19$ ---
see \cite{BH} --- so in this case if $\Gamma(p,q)$ is a Kleinian group
with $p/q=17,18,19$ then it has a subgroup of finite index generated
by rotations.)

For part (iii) of the theorem suppose that $(S^3,K)/G$ is the
unkot. Then the underlying space of $M_0 (p,q)=M(p,q)/G$ is the result
of $(p,dq)$-Dehn-surgery on this unknot with $d=\vert G\vert$.  That
is, $M_0(p,q)$ has underlying space the lens space $L(p,dq)$, which
has a simply-connected $d$-fold covering, proving part (iii). The
proof of (iv) is similar: in this case $M_0(p,q)$ has underlying space
equal to the result of $(p,dq)$ Dehn surgery on a torus knot, which is
a Seifert fibered manifold with infinite fundamental group for
infinitely many $(p,q)$ and with finite fundamental group for
infinitely many $(p,q)$.

Finally, for part (ii), suppose $(S^3,K)$ is invertible.  That is,
there is an involution of $S^3$ which reverses $K$.  The quotient
$S^3/C_2$ by this involution, as a space, is $S^3$, while the quotient
of a tubular neighborhood of $K$ is a ball.  Dehn surgery just
replaces this ball by another ball with different orbifold structure,
so $M(p,q)/C_2$, as a space, is still $S^3$.  The corresponding
$C_2$-extension of $\Gamma$ is thus generated by rotations.\qed

\Section 2.~~Hyperbolic knot complements and groups
generated by rotations

One can also ask when a hyperbolic knot complement $M=S^3-K$ is itself
commensurable with a $\H^3/\Lambda$ with $\Lambda$ generated by
rotations.  Again, this means $\H^3/\Lambda$ has simply connected
underlying space.  If the knot is invertible, then the quotient of
$M=S^3-K$ by the inversion has underlying space an open disk, so the
answer is ``yes.''  Otherwise, $M$ is non-arithmetic by Reid \cite{R}
(who shows that the figure-eight knot is the only knot with arithmetic
hyperbolic complement; the figure-eight knot is invertible) and we can
again argue that the ``orientable commensurator quotient'' $M_0$ of
$M$ (i.e., the quotient of $\H^3$ by the largest Kleinian group
$\Gamma_0$ containing $\Gamma=\pi_1(M)$) would have underlying space
with finite fundamental group.  This could only happen if $\Gamma_0$
is larger than the normalizer $N(\Gamma)$ of $\Gamma$ in
$\PSL(2,\C)$ (since $N(\Gamma)/\Gamma$ is cyclic, by the same
argument as in Sect.~1, so $\H^3/N(\Gamma)$ still has infinite
homology).  As described in \cite{NR, Sect.~9}, this is an exceedingly
rare phenomenon which quite possibly only happens for the figure-eight
knot and the two ``dodecahedral knots'' of Aitcheson and Rubinstein
\cite {AR}.  One of the dodecahedral knots is invertible (this knot is
number $5$ in the series of knots of which number $3$ is shown in
Fig.~1). The other dodecahedral knot is non-invertible and its
$\H^3/\Gamma_0$ is contractible.  Thus summarizing:

\Proclamation Theorem 2.
Let $(S^3,K)$ be a hyperbolic knot.  Then $\Gamma=\pi_1(S^3-K)$ is
commensurable with a Kleinian group generated by rotations if
$(S^3,K)$ is invertible or is the non-invertible dodecahedral knot.
Any other example would have to have the normalizer $N(\Gamma)$ of
$\Gamma$ not equal to the commensurator $\Gamma_0$ and conjecturally
there are no further examples of this.\qed

\Section
3.~~Hidden rotations in Kleinian groups

One can ask whether a group $\Lambda$ commensurable with a given
Kleinian group $\Gamma\subset PSL(2,\C )$ can contain rotations than
$\Gamma$ does not contain --- we call these ``hidden rotations'' for
$\Gamma$.

Of course, Kleinian groups can have hidden rotations --- any
torsion-free subgroup $\Gamma$ of a group $\Lambda$ with torsion does,
for example.  But we can exclude a lot of possible hidden rotations
too.

Let $k(\Gamma )$ be the invariant trace field of $\Gamma$.  That is,
it is the field generated by the traces of squares of elements of
$\Gamma$, cf. \cite{NR}.

\Proclamation Theorem 3.
If a group commensurable with $\Gamma$ contains a $(2\pi
/p)$-rotation, then $\cos (2\pi /p)\in k(\Gamma )$.

\Proof
Suppose $g$ is a $(2\pi /p)$-rotation in a group $\Lambda$
commensurable with $\Gamma$.  Then the trace of $g^2$ is $2\cos (2\pi
/p)$, so $2\cos (2\pi /p)\in k(\Lambda )$.  Since $k(\Lambda)
=k(\Gamma )$ (cf. \cite{NR}), the Theorem follows.\qed

Note that $\cos (2\pi /p)$ is rational for $p\le 4$ and $p=6$, so
Theorem 2 never excludes elements of order $\le 4$ or of order $6$.
However, ``most'' fields will not contain $\cos (2\pi /p)$ for $p=5$
or $p>6$, so ``most'' Kleinian groups $\Gamma$ will admit no hidden
rotations of these orders.

\References{References}

\litemindent= 1 cm  %%%%% amount of indentation for \litem
\ref AR I.~R. Aitchison, and J.~H. Rubinstein, Combinatorial cubings,
cusps and the dodecahedral knots. in {\it Topology 90, Proceedings of
the Research Semester in Low Dimensional Topology at Ohio State}
(Walter de Gruyter Verlag, Berlin - New York 1992), 17--26.

\ref An {E.~M.~Andreev, On convex polyhedra in Lobacevskii space, 
Math. USSR Sbornik 10 (1970), 413--440.}

\ref Ar {M.~A.~Armstrong, The fundamental groups of the orbit space of a 
discontinuous group, Proc. Cambridge Philos. Soc. 64 (1968),
299--301.}

\ref BH S.~A.~Bleiler and C.~D. Hodgson, Spherical space forms and
Dehn fillings, Topology (to appear).

\ref {B} {A. Borel, Commensurability classes and volumes of hyperbolic 
3-manifolds.  Ann. Scuola Norm. Sup. Pisa 8 (1981), 1--33.}

\ref {MB} Morgan,~J.~W. and Bass,~H, {\it eds}, {\it The Smith
Conjecture}, (Academic Press 1984).

\ref {NR} {W. D. Neumann and A. Reid, Arithmetic of hyperbolic
3-manifolds, in {\it Topology 90, Proceedings of the Research Semester
in Low Dimensional Topology at Ohio State} (Walter de Gruyter Verlag,
Berlin - New York 1992), 273--310.}

\ref {NZ} {W. D. Neumann and D. Zagier, Volumes of hyperbolic 3-manifolds, 
Topology 24 (1985), 307--332.}

\ref {R} {A. W. Reid, Arithmeticity of knot complements, J. London
Math. Soc.  (2) 43 (1991), 171--184.}

\ref {T} {W. P. Thurston, The geometry and topology of 3-manifolds,
Mimeographed Notes, Princeton Univ. (1977).}

\ref {V} {E. Vinberg, Discrete groups generated by reflections in
Lobachevskii space, Math. Sb. 114 (1967), 471--488.}

\ref {Z} R. Zimmer, Ergodic Theory and Semi-Simple Lie Groups,
Birkhauser, Boston, 1984.

\end

%% file: degmacs.tex
%%degmacs.tex
%%% Macros for Proceedings Volumes%%%%%%%%%%%%%%%%%%%%%%%%
%%% created June 1992 by Danny Lewis for%%%%%%%%%%%%%%%%%%
%%% Verlag Walter De Gruyter%%%%%%%%%%%%%%%%%%%%%%%%%%%%%%
%%% Revised version July 1994%%%%%%%%%%%%%%%%%%%%%%%%%%%%%

%%%PARAMETERS

\catcode`@=11
\nopagenumbers
\vbadness=10000
\hbadness=10000
\hsize=12.5cm
\vsize=19.2cm
\topskip=12pt
\parindent=0.5cm
\newskip\litemindent
\litemindent= 0,7 cm  %%%%% amount of indentation for \litem
\parskip=0pt
\widowpenalty=10000
\clubpenalty=10000
\hfuzz=1.5pt
\abovedisplayskip=6pt plus 1pt
\abovedisplayshortskip=6pt plus 1pt
\belowdisplayskip=6pt plus 1pt
\belowdisplayshortskip=6pt plus 1pt
\frenchspacing

%%%%FONTS AND FAMILIES%%%%%

% AMS fonts are included, but commented out with
% the symbol "%.%", since not everyone has them.
% If you have them and need them,
% you can remove the %.%.

\font\authorfont=cmti10 at 14.40pt

\font\seventeenrm=cmbx10 at 17.28pt
\font\seventeenit=cmbxti10 at 17.28pt
\font\seventeeni=cmmib10 at 17.28pt
\font\seventeensy=cmbsy10 at 17.28pt
\font\seventeenex=cmex10 at 17.28pt
\font\seventeenmsa=msam10 at 17.28pt
\font\seventeenmsb=msbm10 at 17.28pt
\font\seventeeneuf=eufm10 at 17.28pt

\font\fourteenrm=cmbx10 at 14.40pt
\font\fourteenit=cmbxti10 at 14.40pt
\font\fourteeni=cmmib10 at 14.40pt
\font\fourteensy=cmbsy10 at 14.40pt
\font\fourteenex=cmex10 at 14.40pt
\font\fourteenmsa=msam10 at 14.40pt
\font\fourteenmsb=msbm10 at 14.40pt
\font\fourteeneuf=eufm10 at 14.40pt

\font\twelverm=cmbx10 at 12pt
\font\twelveit=cmbxti10 at 12pt
\font\twelvei=cmmib10 at 12pt
\font\twelvesy=cmbsy10 at 12pt
\font\twelveex=cmex10 at 12pt
\font\twelvemsa=msam10 at 12pt
\font\twelvemsb=msbm10 at 12pt
\font\twelveeuf=eufm10 at 12pt

\font\tenrm=cmr10
\font\tenit=cmti10
\font\tenbf=cmbx10
\font\tenib=cmmib10
\font\teni=cmmi10
\font\tensy=cmsy10
\font\tenbsy=cmbsy10
\font\tenex=cmex10
\font\tenmsa=msam10
\font\tenmsb=msbm10
\font\teneuf=eufm10

\font\ninerm=cmr9
\font\nineit=cmti9
\font\ninebf=cmbx9
\font\ninei=cmmi9
\font\ninesy=cmsy9

\font\ninemsa=msam9
\font\ninemsb=msbm9
\font\nineeuf=eufm9

\font\eightrm=cmr8
\font\eightit=cmti8
\font\eightbf=cmbx8
\font\eighti=cmmi8
\font\eightib=cmmib8
\font\eightsy=cmsy8
\font\eightbsy=cmbsy8

\font\eightmsa=msam8
\font\eightmsb=msbm8
\font\eighteuf=eufm8

\font\sevenrm=cmr7

\font\sevenbf=cmbx7
\font\seveni=cmmi7

\font\sevensy=cmsy7

\font\sevenmsa=msam7
\font\sevenmsb=msbm7
\font\seveneuf=eufm7

\font\fiverm=cmr5
\font\fivebf=cmbx5
\font\fivesy=cmsy5

\font\fivei=cmmi5

\font\fivemsa=msam5
\font\fivemsb=msbm5
\font\fiveeuf=eufm5

\newfam\msafam%Family 8
\newfam\msbfam%Family 9
\newfam\euffam%Family 10

\def\seventeenpoint{\def\rm{\fam0\seventeenrm}%
\textfont0=\seventeenrm\scriptfont0=\fourteenrm
\scriptscriptfont0=\twelverm
\textfont1=\seventeeni\scriptfont1=\fourteeni
\scriptscriptfont1=\twelvei
\def\mit{\fam1\seventeeni}\def\oldstyle{\fam1\seventeeni}%
\textfont2=\seventeensy\scriptfont2=\fourteensy
\scriptscriptfont2=\twelvesy
\def\cal{\fam2\seventeensy}%
\textfont3=\seventeenex\scriptfont3=\fourteenex
\scriptscriptfont3=\twelveex
\def\it{\fam\itfam\seventeenit}%
\textfont\itfam=\seventeenit
\def\bf{\rm}
\def\msa{\fam\msafam\seventeenmsa}%
\textfont\msafam=\seventeenmsa\scriptfont\msafam=\fourteenmsa
\scriptscriptfont\msafam=\twelvemsa
\def\msb{\fam\msbfam\seventeenmsb}%
\textfont\msbfam=\seventeenmsb\scriptfont\msbfam=\fourteenmsb
\scriptscriptfont\msbfam=\twelvemsb
\def\euf{\fam\euffam\seventeeneuf}%
\textfont\euffam=\seventeeneuf\scriptfont\euffam=\fourteeneuf
\scriptscriptfont\euffam=\twelveeuf
\setbox\strutbox=\hbox{\vrule height16pt depth4pt width0pt}%
\baselineskip=20pt\seventeenrm}

\def\fourteenpoint{\def\rm{\fam0\fourteenrm}%
\textfont0=\fourteenrm\scriptfont0=\twelverm
\scriptscriptfont0=\tenbf
\textfont1=\fourteeni\scriptfont1=\twelvei
\scriptscriptfont1=\tenib
\def\mit{\fam1\fourteeni}\def\oldstyle{\fam1\fourteeni}%
\textfont2=\fourteensy\scriptfont2=\twelvesy
\scriptscriptfont2=\tenbsy
\def\cal{\fam2\fourteensy}%
\textfont3=\fourteenex\scriptfont3=\twelveex
\scriptscriptfont3=\tenex
\def\it{\fam\itfam\fourteenit}%
\textfont\itfam=\fourteenit
\def\bf{\rm}
\def\msa{\fam\msafam\fourteenmsa}%
\textfont\msafam=\fourteenmsa\scriptfont\msafam=\twelvemsa
\scriptscriptfont\msafam=\tenmsa
\def\msb{\fam\msbfam\fourteenmsb}%
\textfont\msbfam=\fourteenmsb\scriptfont\msbfam=\twelvemsb
\scriptscriptfont\msbfam=\tenmsb
\def\euf{\fam\euffam\fourteeneuf}%
\textfont\euffam=\fourteeneuf\scriptfont\euffam=\twelveeuf
\scriptscriptfont\euffam=\teneuf
\setbox\strutbox=\hbox{\vrule height13pt depth4pt width0pt}%
\baselineskip=16pt\fourteenrm}

\def\twelvepoint{\def\rm{\fam0\twelverm}%
\textfont0=\twelverm\scriptfont0=\tenbf
\scriptscriptfont0=\eightbf
\textfont1=\twelvei\scriptfont1=\tenib
\scriptscriptfont1=\eightib
\def\mit{\fam1\twelvei}\def\oldstyle{\fam1\twelvei}%
\textfont2=\twelvesy\scriptfont2=\tenbsy
\scriptscriptfont2=\eightbsy
\def\cal{\fam2\twelvesy}%
\textfont3=\twelveex\scriptfont3=\twelveex
\scriptscriptfont3=\twelveex
\def\it{\fam\itfam\twelveit}%
\textfont\itfam=\twelveit
\def\bf{\rm}%
\def\msa{\fam\msafam\twelvemsa}%
\textfont\msafam=\twelvemsa \scriptfont\msafam=\tenmsa
\scriptscriptfont\msafam=\eightmsa%
\def\msb{\fam\msbfam\twelvemsb}%
\textfont\msbfam=\twelvemsb \scriptfont\msbfam=\tenmsb%
\scriptscriptfont\msbfam=\eightmsb%
\def\euf{\fam\euffam\twelveeuf}%
\textfont\euffam=\twelveeuf \scriptfont\euffam=\teneuf%
\scriptscriptfont\euffam=\eighteuf%
\setbox\strutbox=\hbox{\vrule height10pt depth4pt width0pt}%
\baselineskip=14pt\rm}

\def\tenpoint{\def\rm{\fam0\tenrm}%
\textfont0=\tenrm\scriptfont0=\sevenrm
\scriptscriptfont0=\fiverm
\textfont1=\teni\scriptfont1=\seveni
\scriptscriptfont1=\fivei
\def\mit{\fam1\teni}\def\oldstyle{\fam1\teni}%
\textfont2=\tensy\scriptfont2=\sevensy
\scriptscriptfont2=\fivesy
\def\cal{\fam2\tensy}%
\textfont3=\tenex\scriptfont3=\tenex
\scriptscriptfont3=\tenex
\def\it{\fam\itfam\tenit}%
\textfont\itfam=\tenit
\def\bf{\fam\bffam\tenbf}%
\textfont\bffam=\tenbf\scriptfont\bffam=\sevenbf
\scriptscriptfont\bffam=\fivebf
\def\msa{\fam\msafam\tenmsa}%
\textfont\msafam=\tenmsa \scriptfont\msafam=\sevenmsa
\scriptscriptfont\msafam=\fivemsa
\def\msb{\fam\msbfam\tenmsb}%
\textfont\msbfam=\tenmsb \scriptfont\msbfam=\sevenmsb
\scriptscriptfont\msbfam=\fivemsb
\def\euf{\fam\euffam\teneuf}%
\textfont\euffam=\teneuf \scriptfont\euffam=\seveneuf
\scriptscriptfont\euffam=\fiveeuf
\aline=12pt plus 1pt minus 1pt
\halfaline=6pt plus 1pt minus 1pt
\setbox\strutbox=\hbox{\vrule height8.5pt depth3.5pt width0pt}%
\baselineskip=12pt\rm}

\def\ninepoint{\def\rm{\fam0\ninerm}%
\textfont0=\ninerm \scriptfont0=\sevenrm \scriptscriptfont0=\fiverm
\textfont1=\ninei\scriptfont1=\seveni \scriptscriptfont1=\fivei
\def\mit{\fam1\ninei}\def\oldstyle{\fam1\ninei}%
\textfont2=\ninesy \scriptfont2=\sevensy \scriptscriptfont2=\fivesy
\def\cal{\fam2\ninesy}%
\textfont3=\tenex \scriptfont3=\tenex \scriptscriptfont3=\tenex
\def\it{\fam\itfam\nineit}%
\textfont\itfam=\nineit
\def\bf{\fam\bffam\ninebf}%
\textfont\bffam=\ninebf \scriptfont\bffam=\sevenbf
\scriptscriptfont\bffam=\fivebf
\def\msa{\fam\msafam\ninemsa}%
\textfont\msafam=\ninemsa \scriptfont\msafam=\sevenmsa
\scriptscriptfont\msafam=\fivemsa
\def\msb{\fam\msbfam\ninemsb}%
\textfont\msbfam=\ninemsb \scriptfont\msbfam=\sevenmsb
\scriptscriptfont\msbfam=\fivemsb
\def\euf{\fam\euffam\nineeuf}%
\textfont\euffam=\nineeuf \scriptfont\euffam=\seveneuf
\scriptscriptfont\euffam=\fiveeuf%
\aline=11pt plus 1pt minus 1pt
\halfaline=5pt plus 1pt minus 1pt
\setbox\strutbox=\hbox{\vrule height7pt depth3pt width0pt}%
\baselineskip=11pt\rm}

\def\eightpoint{\def\rm{\fam0\eightrm}%
\textfont0=\eightrm \scriptfont0=\sevenrm \scriptscriptfont0=\fiverm
\textfont1=\eighti\scriptfont1=\seveni \scriptscriptfont1=\fivei
\def\mit{\fam1\eighti}\def\oldstyle{\fam1\eighti}%
\textfont2=\eightsy \scriptfont2=\sevensy \scriptscriptfont2=\fivesy
\def\cal{\fam2\eightsy}%
\textfont3=\tenex \scriptfont3=\tenex \scriptscriptfont3=\tenex
\def\it{\fam\itfam\eightit}%
\textfont\itfam=\eightit
\def\bf{\fam\bffam\eightbf}%
\textfont\bffam=\eightbf \scriptfont\bffam=\sevenbf
\scriptscriptfont\bffam=\fivebf
\def\msa{\fam\msafam\eightmsa}%
\textfont\msafam=\eightmsa \scriptfont\msafam=\sevenmsa
\scriptscriptfont\msafam=\fivemsa
\def\msb{\fam\msbfam\eightmsb}%
\textfont\msbfam=\eightmsb \scriptfont\msbfam=\sevenmsb
\scriptscriptfont\msbfam=\fivemsb
\def\euf{\fam\euffam\eighteuf}%
\textfont\euffam=\eighteuf \scriptfont\euffam=\seveneuf
\scriptscriptfont\euffam=\fiveeuf%
\aline=10pt plus 1pt minus 1pt
\halfaline=5pt plus 1pt minus 1pt
\setbox\strutbox=\hbox{\vrule height7pt depth3pt width0pt}%
\baselineskip=10pt\rm}

\skewchar\teni='177
\skewchar\ninei='177
\skewchar\eighti='177
\skewchar\seveni='177
\skewchar\fivei='177
\skewchar\tensy='60
\skewchar\ninesy='60
\skewchar\eightsy='60
\skewchar\sevensy='60
\skewchar\fivesy='60

%%%%headings

\long\def\Title#1\par{%
\global\titlepagetrue
{\parindent=0pt
\raggedcenter\pretolerance=10000
\seventeenpoint #1\par}
\vskip24pt}

\long\def\Author#1\par{%
\centerline{\fourteenpoint\authorfont #1}
\vskip60pt
\tenpoint\noindent\ignorereturn}

\def\Section{\removelastskip
\goodbreak\vskip36pt plus 1pt minus 6pt \section}

\def\section#1\par{%
{\raggedcenter
\interlinepenalty=10000\pretolerance=10000
\noindent\fourteenpoint #1\nobreak\par\nobreak}
\nobreak\vskip12pt\nobreak
\noindent\tenpoint\rm\ignorereturn}

\def\subsection#1\par{%
{\raggedcenter
\interlinepenalty=10000\pretolerance=10000
\noindent\twelvepoint #1\nobreak\par\nobreak}
\nobreak\vskip12pt\nobreak
\noindent\tenpoint\rm\ignorereturn}

\def\References #1 {\ifdim\lastskip<25pt \removelastskip
\vskip24pt plus 1pt\fi
\setbox0=\hbox{\ninepoint [#1]\enspace}%
\centerline{\twelverm References}\nobreak\par\nobreak
\interlinepenalty=10000\parskip=0pt plus 1pt
\litemindent=\wd0
\ninepoint
\nobreak\vskip12pt\nobreak}

\def\runningheads#1#2{\global\outer\def\lefthead{#2}%
\global\outer\def\righthead{#1}}

%%%%GERMAN%%%%%

%%%%FOOTNOTES%%%%%

\newbox\footbox
\setbox\footbox=\hbox{\ninerm 22)~}

\newdimen\footnotespace
\footnotespace=0pt

\def\footnote#1{\let\@sf\empty %#2 (the text) is read later
  \ifhmode\edef\@sf{\spacefactor\the\spacefactor}\/\fi
$^{#1}$\@sf\vfootnote{#1}}

\def\vfootnote#1{\insert\footins\bgroup
\interlinepenalty\interfootnotelinepenalty
\splittopskip\ht\strutbox % top baseline for broken footnotes
\splitmaxdepth\dp\strutbox \floatingpenalty\@MM
\leftskip\z@skip \rightskip\z@skip \spaceskip\z@skip \xspaceskip\z@skip\parindent=\wd\footbox
\litem{\ninerm #1}\ninepoint\footstrut\futurelet\next\fo@t}
\def\fo@t{\ifcat\bgroup\noexpand\next \let\next\f@@t
\else\let\next\f@t\fi \next}
\def\f@@t{\bgroup\aftergroup\@foot\let\next}
\def\f@t#1{#1\@foot}
\def\@foot{\strut\egroup}
\def\footstrut{\vbox to\splittopskip{}}
\skip\footins=\bigskipamount % space added when footnote is present
\count\footins=1000 % footnote magnification factor (1 to 1)
\dimen\footins=8in % maximum footnotes per page

%%%%%OUTPUT ROUTINE%%%%%%

\newif\iftitlepage

\def\makeheadline{\iftitlepage \global\titlepagefalse
\vbox{\titleheadline}%
\vskip12pt
\else \vbox{\ifodd\pageno\rightheadline\else\leftheadline\fi}%
\vskip12pt\fi}

\def\lefthead{}
\def\righthead{}

\def\rightheadline{\line{\vbox to 8.5pt{}%
\ninepoint\hfill\righthead\hfill{\tenrm\folio}}}
\def\leftheadline{\line{\vbox to 8.5pt{}%
\ninepoint{\tenrm\folio}\hfill\lefthead\hfill}}
\def\titleheadline{\line{\vbox to 8.5pt{}%
\ninepoint\it\hfill[Page \folio]}}

\output{\plainoutput}

%%%%%%MACROS%%%%%%%

\def\Litem#1#2{\par\noindent\hangindent#1\litemindent
\hbox to #1\litemindent{\hfill\hbox to \litemindent
{#2 \hfill}}\ignorespaces}
\def\litem{\Litem1}

\newskip\aline \newskip\halfaline
\aline=12pt plus 1pt minus 1pt
\halfaline=6pt plus 1pt minus 1pt
\def\skipaline{\vskip\aline}

\def\ignorereturn{\def\neext{\afterassignment\restart
       \let\next}\neext}
\def\restart{\ifx\next\par\else\let\neext\next\fi\neext}

\def\raggedcenter{\leftskip=0pt plus 4em \rightskip
=\leftskip \parfillskip=0pt \spaceskip=.3333em
\xspaceskip=.5em \pretolerance=9999 \tolerance=9999
\hyphenpenalty=9999 \exhyphenpenalty=9999 }

\def\dotfill{\leaders\hbox to 1em{\hss.\hss}\hfill}

\def\Classification{\ifdim\lastskip<\aline\removelastskip\skipaline\fi
\noindent 1991 Mathematics Subject Classification: }

\def\Abstract{\ninepoint\noindent\bf Abstract. \rm}

\long\def\Proclamation#1. #2\par{\ifdim\lastskip<\aline\removelastskip
\penalty-250 \skipaline\fi{\def\\##1){\litem{\rm(##1)}}\noindent
\bf#1\unskip. \it#2\par}\skipaline}

\def\Proof{\ifdim\lastskip<\aline\removelastskip\skipaline\fi
\noindent\it Proof. \rm}

\def\qedbox{$\rlap{$\sqcap$}\sqcup$}
\def\qed{\nobreak\hfill\penalty250 \hbox{}\nobreak\hfill\qedbox\skipaline}

\def\proclamation#1. {\ifdim\lastskip<\aline\removelastskip\penalty-250
      \skipaline\fi\noindent\bf#1\unskip. \rm}

\def\ref #1 {\litem{[#1]}}

\catcode`@=12
\def\@{\hbox{-}}

%% end of degmacs.tex